\documentclass[12pt]{amsart}
\usepackage{amssymb,amsmath}

\newtheorem{thm}{Theorem}[section]
\newtheorem{lem}{Lemma}[section]

\newtheorem{Example}{Example}

\newtheorem{dfn}{Definition}[section]

\numberwithin{equation}{section}

\def\qed{\hfill $\square$}
\def\comp{\leavevmode\raise.2ex\hbox{${\scriptstyle\mathchar"020E}$}}

\def\P{\mathbf{P}}
\def\E{\mathbf{E}}

\def\proof{\noindent{\it Proof.\ \ }}

\def\FF{\mathcal{F}}
\def\EE{\mathcal{E}}

\begin{document}

\title{Irreducibility and uniqueness of stationary distribution}
\author{Ping He$^1$, Jiangang Ying$^2$}
\address{Department of Applied mathematics, Shanghai University of
Finance and Economics, Shanghai, China} \email{pinghe@shufe.edu.cn}
\address{Institute of Mathematics, Fudan University, Shanghai, China}
 \email{jgying@fudan.edu.cn} \thanks{$^1$supported by
%Research Foundation for Doctor Programme (Grant No. 20040001034) and
National Natural Science Foundation of China (Grant No. 10771131)}
\thanks{$^2$supported in part by NSFC
Grant No. 10671036} \maketitle
\begin{abstract} In this paper, we shall prove that the
irreducibility in the sense of fine topology implies the uniqueness
of invariant probability measures. It is also proven that this
irreducibility is strictly weaker than the strong Feller property
plus irreducibility in the sense of original topology, which is the
usual uniqueness condition.
\end{abstract}

\noindent{\bf 2000 MR subject classification}: 60J45


\section{Introduction}
The existence and uniqueness of invariant measures have been one
of the most important problems in theory of Markov processes. Let
$(P_t)$ be a transition semigroup on a measurable space $(E,\EE)$.
A $\sigma$-finite measure $\mu$ is invariant if $\mu P_t=\mu$ for
any $t>0$, where $$\mu P_t(A):=\int \mu(dx)P_t(x,A).$$ An
invariant probability measure is also called an invariant
distribution or stationary distribution. The existence of an
invariant distribution usually means the positive recurrence and
the uniqueness means ergodicity.

It is well-known (see e.g. \cite{[Doob]}, \cite{[DPZ]},
\cite{[Khas]}) and also very useful that for a nice Markov process
on a nice topological space, the strong Feller property ($P_t$
takes bounded measurable function to continuous function) together
with the irreducibility (any point can reach any open set) implies
the uniqueness of invariant distribution. Usually the
irreducibility is intuitive and not very hard to check. However it
seems that the strong Feller is really strong in many cases
especially in degenerate cases. Besides, two conditions involve
the topology much more than the invariant measure itself does, and
therefore are not so essential.

In this paper we are going to introduce another irreducibility
which is more natural in some sense. For example it depends on the
topology induced by the process itself. We shall prove that the
irreducibility implies the uniqueness of invariant distribution
and also prove that it is really weaker than the strong Feller
property plus the irreducibility (in the sense of original
topology). We also give some characterizations which are easy to
check.

\section{Main results}

Let $$X=(\Omega,\FF,\FF_t,X_t,\theta_t,\P^x)$$ be a right Markov
process on $(E,\EE)$ (say, Polish), with transition semigroup
$(P_t)$ and resolvent (or potential operator) $(U^\alpha:\alpha\ge
0)$. For any nearly Borel set $B$, $T_B$ always denotes the
hitting time of $B$.

\begin{dfn} $X$ is called irreducible if $\P^x(T_G<\infty)>0$ for any $x\in E$ and finely
open $G$, weakly irreducible if $\P^x(T_G<\infty)>0$ for any $x\in
E$ and open $G$.
\end{dfn}

The weak irreducibility is weaker than the irreducibility.

\begin{Example} \rm Let $\nu$ be a probability measure charging on all non-zero
rationals and $X$ the compound Poisson process with L\'evy measure
$\nu$. Then any point is finely open. Since rational numbers are
dense, $X$ is weakly irreducible, but not irreducible with reason
that $X$, staring from 0, can only reach rational numbers.
\end{Example}

\begin{thm} The following are equivalent

$1.$ $X$ is irreducible;

$2.$ For any $A\in\EE$, $U^\alpha1_A$ is either 0 identically or
positive everywhere;

$3.$ All non-trivial excessive measures are equivalent.

\end{thm}

\proof We may assume $\alpha =0$. Suppose (1) is true.  If $U 1_A$
is not identically zero, then there exists $\delta>0$ such that
$D:=\{U 1_A>\delta\}$ is non-empty. Since $U 1_A$ is excessive and
thus finely continuous, $D$ is finely open and the fine closure of
$D$ is contained in $\{U 1_A\ge \delta\}$. Then by Proposition
II.2.8 and Theorem I.11.4\cite{BG},
$$U 1_A(x)\ge P_D U 1_A(x)=\E^x\left( U 1_A(X_{T_D})\right)\ge
\delta \P^x(T_D<\infty)>0.$$ Conversely suppose (2) is true. Then
for any finely open set $D$, by the right continuity of $X$, $U
1_D(x)>0$ for any $x\in D$. Therefore $U 1_D$ is positive
everywhere on $E$.

Let $\xi$ be a non-trivial excessive measure. Since $\alpha \xi
U^\alpha\le \xi$, $\xi(A)=0$ implies that $\xi U^{\alpha}(A)=0$.
However $\xi$ is non-trivial. Thus it follows from (2) that
$U^{\alpha}1_A\equiv 0$, i.e., $A$ is potential zero. Conversely
if $A$ is potential zero, then $\xi(A)=0$ for any excessive
measure $\xi$. Therefore (2) and (3) are equivalent. \qed

It is well-known that the strong Feller property and (weak)
irreducibility together imply the uniqueness of invariant
distribution, which implies the ergodicity. By the strong Feller,
we mean that $P_t$ takes bounded measurable functions to
continuous functions. A condition obviously weaker than strong
Feller is called LSC, which means that for any measurable set $B$,
$U^{\alpha}(\cdot,B)$ is lower semi-continuous. The Brownian
motion is strong Feller, but compound Poisson process is not
strong Feller.

\begin{lem} If $X$ satisfies LSC and weak
irreducibility, then it is irreducible.\end{lem}

\proof Let $A\in\EE$. $U^\alpha 1_A\not=0$ identically. There is
$b>0$ such that $G=\{U^\alpha 1_A>b\}\not=\emptyset$ and is open
due to the property LSC. Since $U^\alpha 1_A$ is
$\alpha$-excessive, we have by Proposition II.2.8 \cite{BG} for
any $x\in E$,
$$U^\alpha 1_A(x)\ge P^{\alpha}_G U^{\alpha}
1_A(x)=\P^x\left(e^{-\alpha T_G}\cdot U^\alpha
1_A(X(T_G))\right).$$ But $X(T_G)\in \bar G$ by Theorem
I.11.4\cite{BG} and then $X(T_G)\ge b$. Hence by the weak
irreducibility, we have
$$U^{\alpha}1_A(x)\ge b \E^x \left(e^{-\alpha
T_G},{T_G<\infty}\right)>0.$$ \qed

\begin{Example} \rm
Let $N=(N_t)$ be a Poisson process with parameter $\lambda>0$ and
$X_t=N_t-t$. Then $X$ does not satisfy strong Feller but satisfies
LSC. Hence it is irreducible. Since $X$ jumps forward and drifts
backward, any set such as $(a,b]$ is finely open. It may be shown
that it satisfies a stronger irreducibility $\P^x(T_y<\infty)>0$ for
any $x,y\in \mathbb{R}$.
\end{Example}

The strong Feller property is certainly too much and hence is not
essential for uniqueness of stationary distribution as indicated in
the following very simple example.

\begin{Example} \rm
Consider the uniform translation $X$ on unit circle. Then $X$ is not
strong Fellerian. However the uniform distribution on circle is the
only stationary distribution of $X$.
\end{Example}

 It seems that in studying uniqueness
problems, our irreducibility is more natural than the irreducibility
under the original topology. In the joint paper \cite{FYZ} of the
second author, it is proven that our irreducibility implies the
uniqueness of ($\sigma$-finite) symmetrizing measures for Markov
processes. Here we shall prove the uniqueness of invariant
distribution under irreducibility.

\begin{thm} The irreducibility of $X$ implies the uniqueness of invariant
distribution.
\end{thm}

\proof We complete the proof in 4 steps.

\begin{itemize}
\item[1.] (\cite{[DPZ]}, Theorem 3.2.4) If a probability measure $\mu$ is invariant, then $\mu$ is ergodic if and only
if for any $A\in\EE$ and $t>0$, $P_t1_A=1_A$, $\mu$-a.s. implies
that $\mu(A)=0$ or 1.

\item[2.] (\cite{[DPZ]}, Theorem 3.2.5) If $\mu$ and $\nu$ are ergodic to
$(P_t)$ and $\mu\not=\nu$, then $\mu$ and $\nu$ are singular.

Note that the condition for 1 and 2 is the stochastic continuity
of $X$ which is surely implied by the right continuity of $X$.

\item[3.] If $(P_t)$ is irreducible, then any invariant
distribution is ergodic.

In fact for any $A\in\EE$, $t>0$, $$P_t1_A=1_A,\ \mu-\text{a.s.}$$
Then for any $B\in \EE$, $$\langle 1_B,P_t1_A\rangle_\mu=\langle
1_B,1_A\rangle_\mu$$ for any $t>0$. It follows by Fubini's theorem
that
$$\langle 1_B,\alpha U^\alpha 1_A\rangle_\mu=\langle
1_B,1_A\rangle_\mu$$ for any $\alpha>0$. Hence $\alpha U^\alpha
1_A=1_A$, $\mu$-a.s. Our irreducibility is equivalent to $U^\alpha
1_A$ is either 0 identically or positive everywhere. This implies
that $\mu(A)=0$ or 1.

\item[4.] Irreducibility implies the uniqueness of invariant
distribution.

In fact, if $\mu$ and $\nu$ are two different invariant
distributions, then they are ergodic by (3) and they are singular
by (2). However the irreducibility implies that all excessive
measures are equivalent to each other, which leads to a
contradiction.

\end{itemize}

That completes the proof. \qed

\noindent{\bf Acknowledgement:} The authors would like to thank
Professor Dong, Zhao of Academia Sinica for stimulating discussion.

\end{document}